\documentclass{article}
\usepackage{amsmath,amsthm,amsfonts,amscd,amssymb,eucal,latexsym,hyperref} 
\usepackage[all]{xy}
\begin{document}

\title{A rigidity result for  extensions of braided
 tensor  $C^*$--categories derived from compact matrix quantum groups}
  
\author{Claudia Pinzari  and John E. Roberts}
\date{}
\maketitle

   \centerline{\it Dedicated to Sergio Doplicher on the occasion of his seventieth birthday}                
\bigskip

\begin{abstract}
Let $G$ be a classical compact Lie group and $G_\mu$ the associated compact matrix quantum group  deformed by a  positive parameter $\mu$ (or $\mu\in{\mathbb R}\setminus\{0\}$ in the type $A$ case).  It is well known that the category of unitary 
representations of $G_\mu$ is a braided tensor $C^*$--category. We show that any braided tensor $^*$--functor $\rho: \text{Rep}(G_\mu)\to{\cal M}$ to another braided tensor $C^*$--category with irreducible tensor unit   is full if $|\mu|\neq1$. In particular, the functor of restriction  Rep$G_\mu\to\rm{Rep}(K)$ to a proper compact quantum subgroup $K$ cannot be made into a braided functor.
Our result also shows that  the Temperley--Lieb category ${\cal T}_{\pm d}$ for $d>2$ can not be embedded properly into a larger category with the same objects as a braided tensor $C^*$--subcategory.

\end{abstract}

\section{Introduction} 

There are various  strategies for studying the structure of or classifying semisimple rigid  
tensor categories. As an oversimplification,
sometimes one focuses on the structure of simple objects.
A basic result is Ocneanu's 
rigidity, see \cite{ENO} for a proof, asserting that there  are finitely 
many ${\mathbb C}$--linear fusion categories with 
a prescribed fusion ring. Important results have been obtained for  categories that may or may not have finitely many irreducibles, or a braiding, see   \cite{KW, Muger, O1, O2, RSW, TW, ENO}, however this is  an incomplete list. 

Another approach, when there is
a relevant subcategory with few arrows, is to try and construct 
the whole category ${\cal M}$ from a 
smaller subcategory ${\cal A}$.
This subcategory is regarded as a {\it symmetry} of ${\cal M}$, and the classification problem becomes that of classifying extensions ${\cal M}$ with the given symmetry.
 
The first example, motivated by AQFT \cite{DHR} (see also \cite{Haag}), is that of permutation symmetry. For tensor $C^*$--categories  with conjugates (and subobjects, direct sums and irreducible tensor unit)
realizations of this symmetry are few, they are
classified by a single integer parameter and a sign, the {\it statistics phase}. When all the statistics phases are one, ${\cal M}\simeq\text{Rep}(G)$ with its natural permutation symmetry, for a unique compact group $G$ \cite{DR}. A related independent result is also well known \cite{D}. 

The theory of subfactors \cite{Jones}  or low dim AQFT  \cite{FRS}, although  differing, provide further remarkable instances of this general scheme. In the first case the Jones projections play a fundamental role, leading 
to  the Temperley--Lieb symmetry; whereas  the  Virasoro algebra symmetry  with central charge $<1$ plays a major role in  CFT on the circle.
Deep classification results have been obtained in both these areas \cite{GHJ, H,  Popa, AH},   \cite{KL}.

For the categories arising from low dimensional 
QFT or from certain quantum groups one has braid group 
symmetries. The main difference from the permutation symmetric case is the great variety
of realizations in tensor categories.

In this paper we start with braided categories ${\cal A}$ arising from quantum groups where the braiding comes from an $R$--matrix.
The matrix carries more information on the quantum group than in the group case, where
most of the  information gets lost, so that, for example, for any proper closed subgroup $K$ of a group $G$ the restriction functor $\text{Rep}(G)\to\text{Rep}(K)$ is a not full
permutation symmetric tensor $^*$--functor 
$\text{Rep}(G)\to\text{Rep}(K)$.

The $R$--matrix  depends on  the structural constants of the quantum group,  
making it an object intimately related to the quantum group and
raising the question of whether
the class of  tensor categories with that specific braided symmetry  has a 
rigid structure. 

A recent result in \cite{CKL} for the categories 
arising from CFT  on the circle seems  to confirm a certain rigidity when constructing 
inclusions of unitarily braided, but not permutation symmetric, 
tensor categories.

More precisely, the authors, working with 
inclusions of conformal nets on the circle,
where  unitary braided symmetries arise, have shown that
far fewer low values of the Jones index of the inclusion can occur than in 
subfactor situation. They exclude all 
non-integer index 
values $<3+\sqrt{3}$ (hence in particular all non-integer index values in the Jones 
discrete series) 
except one,
$4\cos^2\pi/10$, which is realized \cite{Rehren}.
The integer values are known to be realized from inclusions derived by taking 
fixed points under  group actions.

$\alpha$--induction relates the tensor category associated with a net
${\cal A}$ to that associated with an extension ${\cal B}$, see \cite{KL} and references there.  
We also mention the work of \cite{muger} for a categorical  relation corresponding  
to inclusions of nets.

The aim of this paper is to derive a  rigidity result for 
 tensor categories with a braided symmetry derived from 
quantum groups of all Lie types but not at a 
root of unity. Note that  in this  case the $R$--matrix  carries
maximal information on the quantum group.

We shall work with tensor $C^*$--categories,   ${\cal 
A}$ will 
be the category of unitary representations of 
any  deformation $G_\mu$ of the classical compact 
Lie groups $G$ by a  real  parameter $\mu$.

It is known that  the  $R$--matrix for such quantum groups
makes the associated representation category into a braided tensor $C^*$--category
with conjugates, the analogues of duals. Except for the extreme values of $\mu$, the coinverse of $G_\mu$  
is not involutive and the braiding not unitary and we show that every braided tensor $^*$--functor 
$\rho:\text{Rep}(G_\mu)\to {\cal M}$ to a braided tensor 
$C^*$--category is full (see Theorem 3.1 for a more general statement). 
We derive this result from Theorem 5.4, showing an obstruction to construct extensions 
of braided tensor $C^*$--categories.

Note that, in our class of examples, ${\cal A}$ has  infinitely many 
irreducibles whose indices (or dimensions) can be arbitrarily large.

The category ${\cal M}$ is not assumed to be embeddable into the category of Hilbert spaces.
However, when it is, our result shows that  the braiding of Rep$(G_\mu)$ does not extend to any
Rep$K$, where $K$ is a proper quantum subgroup.

This corollary does not, to our knowledge, seem to have been noticed  in the 
literature. It shows how the results are discontinuous in the classical limit  $\mu= 1$.

Another application concerns Temperley--Lieb categories. It states that when this category is generated by an object of intrinsic dimension $d>2$,
(i.e. not in the discrete series of Jones), it cannot be 
embedded properly as a braided tensor $C^*$--subcategory into a larger such category with the same objects.
We recall here the situation for subfactors 
$N\subset M$ with finite index where the associated $N$--$N$--bimodule tensor category 
contains a copy of the TL category  as a (non-braided) tensor 
$C^*$--category.

The relation to our result can be seen by regarding the 
 TL category ${\cal T}_{\pm d}$ as a full braided tensor subcategory of the 
representation category of $S_{\mu} U(2)$ for a deformation parameter
determined by $|\mu+\mu^{-1}|=d$, and $\mu$ is either positive or negative corresponding to the minus or plus sign.  
 
Our result relies on the theory of induction of 
 \cite {PR}.  We showed there that if a tensor $C^*$--category with conjugates ${\cal A}$
 admits a tensorial embedding $\tau$ into the category of Hilbert 
spaces, and hence, by Woronowicz duality, 
is a full tensor $C^*$--subcategory of the representation category of a compact quantum group $G^\tau$, 
then  given an embedding of 
$\rho:{\cal A}\to{\cal M}$ of tensor $C^*$--categories,
  the full subcategory of 
${\cal M}$, whose objects are in the image of $\rho$, can be identified 
with  a category of 
Hilbert bimodule representations of $G^\tau$ for
a noncommutative $G^\tau$--ergodic $C^*$--algebra $({\cal C},\alpha)$
intrinsically associated with $\rho$ and $\tau$. $({\cal C},\alpha)$ is to be regarded  
as a virtual subgroup. For a fixed $\tau$, 
isomorphic pairs $({\cal M}, \rho)$,  $({\cal M}', \rho')$ yield conjugate ergodic actions.

We show that, taking braided symmetries into account,
there is an obstruction to constructing such non-trivial ergodic actions $({\cal C},\alpha)$. More precisely, our main result states that   the invariant 
$\kappa$ of ${\cal A}$  must be a phase on irreducible spectral representations 
for the action $\alpha$ corresponding to $\rho:{\cal A}\to{\cal M}$ (Theorem 5.4).
We are thus left to check that $\kappa(v)$ is never a phase for our quantum group $G_\mu$, unless $v$ is equivalent to the trivial representation (Theorem 3.1).

The  invariant $\kappa$  appeared in its simplest form for categories derived from  AQFT.   It is just  the statistics phase 
recalled above for permutation symmetric tensor $C^*$--categories \cite{DHR}. In two spacetime dimensions the braid group replaces the permutation group.
Unitarity of the braiding  means that $\kappa$ is a phase  
 on irreducible objects but it need not be a sign \cite{FRS}. For general braided tensor $C^*$--categories, $\kappa$ has been discussed  in \cite{LR}.

A closely related invariant was independently  
introduced   for  certain braided tensor categories, the ribbon categories \cite{RT}, and referred to as the {\it twist}, $\theta$. 
Its importance derives from the fact that ribbon categories associated with certain quantum groups at
roots of unity lead to the Reshetikhin--Turaev invariants of $3$--manifolds and  links in $3$--manifolds.

Note that for quantum groups  with unitary braided symmetries, 
$\kappa$ is indeed always a phase, and our obstruction vanishes. 
In this case, braided extensions can be expected.

We describe the braided symmetries of $\text{Rep}(S_\mu 
U(d))$ 
got by  twisting  the universal $R$-matrix 
by $d^{\text{th}}$ roots of unity. The 
 invariant $\kappa$ can be used to show that these braidings provide inequivalent 
braided tensor categories.

We finally mention that reconsidering 
the construction of these braidings, led us to compare explicitly the intrinsic characterizations of 
the quantum $SU(d)$, when not at roots of unity,  given by \cite{KW} and \cite{P}.
Whilst \cite{KW} starts from the fusion rules, \cite{P} relies on the braiding. The twist
 $\tau_{\cal C}$ may be derived following the approach of \cite{P}, leading to a characterization of the tensor
 categories with fusion rules of type   $sl(d)$   in the spirit of 
\cite{KW}, in terms of generators and relations analogous to \cite{P}, but 
now involving a $d^{\text{th}}$ root of unity, see Prop.\ 7.2 and Cor.\ 7.3.

The paper is organized as follows. In Sect.\ 2 we recall the notion of a braided 
tensor $C^*$--category.
We shall propose a variant weaker than the usual one, but sufficient for our main result, announced in Sect.\ 3.  To this end, we recall
the natural structure of $\text{Rep}G_\mu$ as a braided tensor $C^*$--category.
We shall in particular emphasize the case of $S_\mu U(d)$ as it allows one to grasp 
the general  idea of the proof quickly. 
 In Sect.\ 4 we derive the main properties of the invariant $\kappa$ of a
 braided tensor $C^*$--category with conjugates and compare it with the twist of a ribbon tensor
 category. In Sect.\ 5 we show that  when the invariant $\kappa$ of a 
braided tensor $C^*$--category ${\cal A}$ is not unitary it 
is an obstruction to construct braided extensions. In Sect.\ 6 we complete 
the proof by showing that  $\kappa(v)$ is never a phase for an 
irreducible representation $v$ of the quantum groups $G_\mu$,  
unless $v$ is the trivial representation. In the last section 
we reconsider the results of \cite{KW} and \cite{P}.

\section{Braided tensor $C^*$--categories} 
\medskip

In this section we introduce  various notions of braiding, from the weakest to the strongest.
These will allow us to isolate that aspect to our main result.
  
 The categories we shall consider   will always be assumed to be strict, to be over the complex numbers, with 
irreducible tensor unit $\iota$. We shall also assume existence of subobjects and direct sums, unless otherwise stated, but see the remarks at the beginning of Sect.\ 4. In this paper we shall mostly deal with tensor $C^*$--categories (or unitary categories).

The weakest  notion   requires that
for any pair of objects $u$, 
$v$ of a tensor $C^*$--category ${\cal A}$ there is an invertible intertwiner 
$\sigma(u,v)\in(u\otimes v, v\otimes u)$ 
such that for arrows $T$ belonging to spaces of the form
$(u,\iota)$,  
 $$1_v\otimes T\circ\sigma(u,v)=
T\otimes 1_{v}.\eqno(2.1)$$ 
Note that the dual intertwiners 
$\sigma_d(u,v):=
 {\sigma(u,v)^*}^{-1}\in(u\otimes v,v\otimes u)$ 
do not satisfy $(2.1)$. 
If both $\sigma$  and $\sigma_d$ satisfy $(2.1)$, $\sigma$ will be referred to as   a {\it weak left   braided symmetry}.
Notice  that a  weak left  braided symmetry by itself,  may  not even be related to the 
braid group 
as we are not assuming  properties $(2.3)$ and 
$(2.4)$ below. However, the relation to the braid group will follow automatically for the weak braided symmetries of interest in this paper (cf.\ the end of the section).

 A stronger  notion is that of a {\it left   braided 
symmetry}.
One requires the
following relations for objects $u$, $u'$, $v$, 
and arrows $T\in(u,u')$, 
 $$\sigma(\iota,u)=\sigma(u,\iota)=1_u,\eqno(2.2)$$
 $$\sigma(u\otimes u', v)=\sigma(u,v)\otimes 
1_{u'}\circ1_{u}\otimes\sigma(u',v),\eqno(2.3)$$
$$\sigma(u',v)\circ T\otimes 1_v=1_v\otimes 
T\circ\sigma(u,v).\eqno(2.4)$$
Indeed, the notion of a weak left  braided symmetry
is just a special case of the naturality property $(2.4)$ and $(2.2)$  
for $u$ or $u'$ the trivial object $\iota$.

A {\it weak right    braided symmetry} is defined by the following equation for $T\in(u,\iota)$,
$$1_v\otimes 
T=T\otimes1_v\circ\sigma(v,u)=T\otimes1_v\circ\sigma_d(v,u),$$
and  {\it right braided symmetry} is defined replacing $(2.3)$ and $(2.4)$ by
$$\sigma(u, 
v\otimes v')=1_v\otimes\sigma(u,v')\circ
\sigma(u,v)\otimes 
1_{v'},\eqno(2.3)'$$
$$\sigma(v,u')\circ1_v\otimes 
T=T\otimes1_v\circ\sigma(v,u).\eqno(2.4)'$$
Obviously, $\sigma$ is a  (weak) right  braided symmetry if and only if
$\sigma_{-1}(u,v):=\sigma(v,u)^{-1}$ 
or $\sigma_*(u,v):=\sigma(v,u)^*$
are (weak)  left braided symmetries.

We thus recover the usual notion of a braided symmetry considering 
 a   left and right braided symmetry.
A  (weak)  half braided symmetry will just mean a  (weak) left or   right braided symmetry.
Given (weak)  half  braided symmetry $\sigma$ in ${\cal A}$, the dual symmetry
$\sigma_d(u,v)$ is another  braided 
symmetry 
for ${\cal A}$ of the same type.

 The representation categories of  $S_\mu U(d)$,  or, more generally, of $G_\mu$, where $G$ is a classical compact Lie group, or of $A_o(F)$ are 
well known examples of braided tensor 
$C^*$--categories. These will in fact be our main examples. They will be discussed in later sections. 

If $\sigma$ satisfies $(2.1)$ or is a weak half  or half braided symmetry for a tensor 
$C^*$--category ${\cal A}$ and 
$\tau:{\cal A}\to\text{Hilb}$ is a tensor $^*$--functor into the category of Hilbert spaces then the representation category of the  
compact quantum group $G^{\tau}$ associated to $\tau$ via Woronowicz duality has a braiding $\sigma_\tau$ of the same type defined by  $\sigma_\tau(\hat{u},\hat{v}):=\tau(\sigma(u,v)),$
where $\hat{u}$ is the representation  of $G^{\tau}$ corresponding to the object $u$ of ${\cal A}$.

It is easy to show  that the category 
  of  Hilbert spaces has a unique braiding, 
even in the sense of $(2.1)$, given by its unique permutation symmetry. 
  If $({\cal A},\sigma)$ and $({\cal M},\varepsilon)$ are tensor $C^*$--categories with some type of braiding, a tensor $^*$--functor $F:{\cal A}\to{\cal M}$ will be called braided
 if for any pair of objects $u$, $v$ in ${\cal A}$, $$F(\sigma(u,v))=\varepsilon(F(u), F(v)).$$
 We shall consider cases where ${\cal A}$ has a stronger braided symmetry than ${\cal M}$.  Note that, if  for example properties $(2.2)$, $(2.3)$ hold for the braiding of ${\cal A}$ 
 and $F$ is a braided tensor ${}^*$--functor to a weak half  braided ${\cal M}$, then those properties also
hold for the braided symmetry of ${\cal M}$ for objects in the image of
 $F$. However, the naturality axioms $(2.4)$ or $(2.4)'$ are not inherited by the full subcategory of ${\cal M}$ whose objects are in the image of $F$.\medskip

 \section{The rigidity result }
 
Tensor categories derived from quantum groups arising from a deformation of classical Lie groups admit  natural braided symmetries associated with 
the universal $R$--matrices of Drinfeld,  see, e.g. \cite{CP, KS}. We work with compact matrix quantum groups. We start by describing briefly the corresponding braided symmetry for the type $A_{d-1}$ in this framework. 
After this, we recall how to construct braided tensor $C^*$--categories for the other Lie types, before illustrating our main result and its corollaries.

Woronowicz has introduced the compact quantum group $S_\mu U(d)$, for a nonzero $\mu\in[-1, 1]$, using his duality theorem \cite{W}. It is the (maximal) compact quantum group whose representation category $\text{Rep}(S_\mu U(d))$ is generated, as an embedded tensor $C^*$--category, by the deformed
determinant element $$S=\sum_{p\in{\mathbb 
P}_d}(-\mu)^{i(p)}\psi_{p(1)}\otimes\dots\otimes\psi_{p(d)},$$ 
where $\psi_i$ is an orthonormal basis of a $d$--dimensional Hilbert space.

For a nonzero complex parameter $q$, let $H_n(q)$ denote the Hecke algebra of type $A_{n-1}$, i.e.
 the quotient of the complex group algebra of the 
braid group
${\mathbb B}_n$ with   generators $g_1,\dots g_{n-1}$ by the
 relations 
 $$g_i^2=(1-q)g_i+q.$$ For $q$ real or $|q|=1$, $H_n(q)$ becomes a $^*$--algebra 
 with involution making the spectral idempotents of the $g_i$ into selfadjoint projections.
 However, there are nontrivial Hilbert space representations for $n$ 
arbitrarily large, only if $q>0$ or $q=e^{\pm2\pi i/\ell}$ \cite{Wenzl}. In 
this paper we only consider the real case.

There is a   well known remarkable representation of the Hecke algebra  $H_\infty(q)$, with $q=\mu^2$,
in $\text{Rep}(S_\mu U(d))$, due to Jimbo and Woronowicz \cite{J,W},
and defined by
 $$\eta(g_1)\psi_i\otimes\psi_j=\mu\psi_j\otimes\psi_i, \quad 
i<j,$$
$$\eta(g_1)\psi_i\otimes\psi_i=\psi_i\otimes\psi_i, $$
$$\eta(g_1)\psi_i\otimes\psi_j=\mu\psi_j\otimes\psi_i+(1-q)\psi_i\otimes\psi_j, 
\quad i>j.$$
$\text{Rep}(S_\mu U(d))$ admits an intrinsic characterization  \cite{KW}, \cite{P}.
  We follow the approach of
 \cite{P} (based on the type of the braiding) and compare with that of \cite{KW} (based on the type of the fusion rules) in the last section.
 
$\text{Rep}(S_\mu U(d))$ is, up to equivalence, the unique tensor $C^*$--category ${\cal A}$ (with subobjects and direct sums) with a tensor 
$^*$--functor from the braid category $\eta:{\mathbb B}\to{\cal A}$ 
factoring through representations of the complex Hecke algebras of type $A$,
$H_n(q)$ for $q=\mu^2$ 
 and generated by
an object $u$ and an arrow $S\in(\iota, u^d)$  satisfying
$$S^*\circ S=d!_q,\quad S\circ S^*=\eta(A_d),\eqno(3.1)$$
$$S^*\otimes 1_u\circ 1_u\otimes S=(d-1)!_q(-\mu)^{d-1},\eqno(3.2)$$
$$\eta(g_1\dots g_d)\circ S\otimes1_u=\mu^{d-1}1_u\otimes 
S,\eqno(3.3)$$
where $n!_q$ is the usual quantum factorial and $A_d$ is the antisymmetrized sum of the 
elements of the canonical basis of $H_d(q)$,  
a scalar multiple of  the analogue of the 
totally antisymmetric projection.  
The above relations are realized by the deformed
determinant element $S$ and the  
JW representation $\eta$. 
They are easily verified for $d=2$.
(see \cite{P} for details. Notice our $g_i$ corresponds 
to $-g_i$ there.)

\medskip

\noindent{\it Remark}
Note that for $d$ odd, $\text{Rep}(S_\mu U(d))$ and $\text{Rep}(S_{-\mu}U(d))$ are canonically isomorphic. \medskip

It is well known from Drinfeld's theory of universal 
$R$--matrices, that the map 
$\sigma_\omega: g_i\to\frac{\omega}{\mu}\eta(g_i)$, where 
$\omega$ is a complex $d^{\text{th}}$ root of $\mu$, makes $\text{Rep}(S_\mu 
U(d))$ into a braided tensor category. It is easy to check that 
$\text{Rep}(S_\mu U(d))$ actually becomes a braided tensor $C^*$--category 
in this way, cf.\ \cite{P}.

Working with a fixed deformation parameter leads to $d$
inequivalent braided symmetries obtained by varying $\omega$.
This fact was first noted in  \cite{KW}, see also \cite{TW}, Sect.\ 4.

\medskip

Note that, for $\mu>0$, the braided symmetry of $\text{Rep}(S_\mu U(d))$
corresponding to the  positive $d^{\text{th}}$ root $\mu^{1/d}$ of $\mu$ is a natural choice, as it reduces
to the unique permutation symmetry of the category of Hilbert spaces for $\mu=1$. \medskip

Deformation has been generalized to all classical compact Lie groups 
 in the framework of compact matrix quantum groups
\cite{R, A, KS}. The starting point was the 
 dual picture of Drinfeld and Jimbo.

  We start with a complex simple Lie algebra $g$.  Denote by $(\alpha_i)_1^r$ a set of simple roots.  Let $A=(\frac{2(\alpha_i,\alpha_j)}{(\alpha_j,\alpha_j)})$ be the Cartan matrix, where $(\cdot{}\ , \cdot{}\ )$  is a symmetric invariant bilinear form on $g$ such that $(\alpha,\alpha)=2$ for a short root $\alpha$. It follows that
$d_j:=\frac{(\alpha_j,\alpha_j)}{2}\in{\mathbb N}$ for all $j$.

For  a complex deformation parameter $\mu$ other than a root of unity,
consider  
 the quantized universal enveloping algebra $U_\mu(g)$ as defined
  in Ch. 2, 7.1.1, \cite{KS}.

The category ${\cal C}(g,\mu)$ of finite 
 dimensional representations of $U_\mu(g)$ admitting a weight decomposition has
 the structure of a 
 braided ribbon tensor category \cite{KS}, see also  \cite{CP}, \cite{BK}.
 
 If $\mu\in(0,1)$, $U_\mu(g)$ becomes a Hopf $^*$--algebra with real form
  $E_i^*=F_i$, $K_i^*=K_i$.
The representation space of every  object of  ${\cal C}(g,\mu)$  with {\it positive} weights
   on the $K_i$'s  admits a natural Hilbert space structure
 making it into a $^*$--representation. One gets in this  
 way a tensor $C^*$--category with conjugates and braided symmetry embedded into the category of Hilbert spaces.  A compact quantum group   $G_\mu$ may then be defined via Woronowicz duality, see \cite{R, A, KS}.
 
 The  category $\text{Rep}(G_\mu)$ of unitary finite dimensional representations of $G_\mu$ will be regarded as a braided tensor $C^*$--category
 as described above. In the particular case of  type $A$, we allow  negative parameters for $S_\mu U(d)$ and the braided symmetry may be any of the $\sigma_\omega$'s.
  
The aim of this paper is to prove the following rigidity result  (a consequence 
of Theorem 5.4.)

\medskip

\noindent{\bf 3.1. Theorem} {\sl  For $\mu\neq1$ (and also $\mu\neq -1$  in the type $A$ case), every braided tensor
$^*$--functor  $\rho:\text{Rep}(G_\mu)\to{\cal M}$ into a weak half  braided tensor $C^*$--category is 
full.}\medskip

We derive the following consequences.
 \medskip

\noindent{\bf 3.2. Corollary} {\sl For $\mu\neq 1$ (and also $\neq -1$ in the type $A$ case) the braided 
symmetry of $\text{Rep}(G_\mu)$ does not make 
 the representation category of any proper compact 
quantum 
subgroup into 
a  weak half  braided tensor $C^*$--category.}\medskip

\noindent{\it Proof} If $K$ is a proper compact quantum subgroup  then the 
 tensor $^*$--functor $\text{Rep}(G_\mu)\to\text{Rep}(K)$ of restriction is not full, hence 
the image of the braided symmetry of $\text{Rep}(G_\mu)$ can not be a braided 
symmetry for $\text{Rep}(K)$.
\medskip

The next corollary is an application to the Temperley--Lieb category ${\cal T}_{\pm d}$.
Recall that for $d>0$, ${\cal T}_{\pm d}$ may be defined as the universal tensor $^*$--category
with objects ${\mathbb N}_0$ whose
 arrows are generated by a single 
arrow $R \in (0, 2 )$ satisfying 
$R^*\otimes 1_1\circ 1_1\otimes R=\pm1_1$ and $R^*\circ R=d$.
The 
following assertions are well known \cite{GW}. The units and generating objects of these 
categories are irreducible and the spaces of arrows are finite dimensional. The 
categories are simple except at roots of unity, $d = 2 \cos\pi/n$, 
 when 
they have a single non-zero proper ideal. They are tensor $C^*$--categories 
when $d \geq 2$ and at roots of unity their quotients by the unique non-zero proper 
ideal are tensor $C^*$--categories. The projection $e=\frac{R\circ R^*}{d}$ defines
a   Hecke algebra representation  $\eta(g_1):=(1+q)e-q$ where $q=\mu^2$ and 
$-1<\mu<1$ is determined by $|\mu+\frac{1}{\mu}|=d$ and $\mu>0$ for ${\cal T}_{-d}$ and $\mu<0$ for ${\cal T}_d$. It is also known that ${\cal T}_{\pm d}$ is a full braided tensor $C^*$--subcategory of $\text{Rep}(S_{\mu}U(2))$, see \cite{PR1} with previous results in   \cite{B} for embedded categories. It is equivalent to $\text{Rep}(S_{\mu}U(2)$ after completion under subobjects and direct sums. We may thus apply Theorem 3.1.

\medskip

\noindent{\bf 3.3. Corollary} {\sl The  category ${\cal T}_{\pm d}$ generated by a 
single selfadjoint object of dimension $d>2$ regarded as a 
braided 
tensor $C^*$--category can not be embedded properly 
as a weak half  braided tensor $C^*$--subcategory with the same objects.}\bigskip

\bigskip

\section{ Invariant $\kappa$ and the twist of ribbon categories }

In this section we recall  the  main properties of the  invariant $\kappa$ of a braided tensor $C^*$--category with conjugates  \cite{LR}   and we compare it 
with the twist $\theta$ of ribbon categories  \cite{RT}.

 If $v$ is an irreducible object of ${\cal A}$ and $R\in(\iota,\overline{v}\otimes v)$, $\overline{R}\in(\iota, v\otimes\overline{v})$ a solution of 
 the conjugate equations for $v$ then $\sigma(v, \overline{v})\circ \overline{R}$ must be a nonzero scalar multiple of $R$ whenever $\sigma(v,\overline{v})\in(v\otimes\overline{v},\overline{v}\otimes v)$ is an invertible intertwiner, as the space $(\iota, \overline{v}\otimes v)$ is one dimensional. 
This scalar  in general depends on the choice of the 
conjugate equations. However, the following proposition shows that it is independent if we choose a {\it standard} solution
of the conjugate equations, in the sense of 
 \cite{LR}.
 For irreducible objects $v$ solutions are standard if and only if they are normalized, $\|R\|=\|\overline R\|$. \medskip

\noindent{\bf 4.1. Proposition} {\sl If $\sigma$ is a right braided symmetry for ${\cal A}$, then for any irreducible object $v$ the scalar $\kappa_r(v)$ defined by $$\sigma(v,\overline{v})\circ\overline{R}=\kappa_r(v) R$$
does not depend on the choice of  
standard solutions  of the conjugate equations for $v$.

 If $\sigma$ is a 
braided symmetry, $\kappa_r$ is a class function, 
i.e. $\kappa_r(v)=\kappa_r(u)$ 
if $v$ and $u$ are unitarily equivalent.}\medskip

\noindent{\it Remark} Note that, unlike in the case of tensor categories,  we do not need to make a choice $v\to R_v$ to define $v\to\kappa_r(v)$. 
\medskip

We may extend
 $\kappa_r$ to reducible objects $v$ of ${\cal A}$ 
by the equation
$$\sigma(v,\overline{v})\circ\overline{R}=1_{\overline{v}}\otimes\kappa_r(v)\circ R,$$
where we again use standard solutions.  Note that $\kappa_r(v)\in(v,v)$.
 Similarly, for left braided symmetries we may define 
$\kappa_l(v)$.

For a braided symmetry,
$\kappa_r(v)=\kappa_l(\overline{v})=:\kappa(v).$
 In this case $\kappa(v)=\sum\kappa(v_n)E_n$ where $v_n$ are the 
irreducible components of $v$ and $E_n$ minimal central
projections in $(v,v)$ with $\sum_nE_n=1_v$. It follows that $v\to\kappa(v)$ is central,
$$\kappa(v)\circ T=T\circ\kappa(u),\quad T\in(u,v),$$
see \cite{LR}.

\medskip

\noindent{\bf 4.2. Proposition} {\sl If $\sigma$ is a right braided symmetry, then
for any object $v$ of ${\cal A}$
$$\kappa^{\sigma_{-1}}_l(v)=\kappa^\sigma_r(v)^{-1},
\quad 
\kappa^{\sigma_*}_l(v)={\kappa^\sigma_r(v)}^*,\quad 
\kappa^{\sigma_d}_r(v)={\kappa^\sigma_r(v)^{-1}}^*.$$}\medskip

\noindent{\it Proof} We may assume $v$ irreducible. The first relation is 
 obvious and the third follows from the first two. Now 
 $\sigma(v,\overline{v})^*\circ 
R_v=
\kappa^{\sigma_*}_l(v) \overline{R}_v$, so
$$\kappa^{\sigma_*}_l(v)\|\overline{R}_v\|^2=
\overline{R}_v^*\circ\sigma(v,\overline{v})^*\circ 
R_v=(\sigma(v,\overline{v})\circ \overline{R}_v)^*\circ 
R_v=\overline{\kappa_r(v)}\|R_v\|^2,$$
completing the proof.
\medskip

The central element $\kappa$
is in fact an invariant for braided tensor $C^*$--categories.
More precisely, the invariant $\kappa_r$ is preserved under a 
full tensor $^*$--functor of right braided tensor $C^*$--categories with
conjugates.

We next give an alternative definition of $\kappa_l(v)$,    in the spirit of the 
 relation between statistics parameter, dimension and statistics phase in AQFT \cite{Haag}.
 \medskip
 
 \noindent{\bf 4.3. Proposition} {\sl If $\sigma$ is a left braided symmetry of ${\cal A}$, then for any 
 object $v$ and any
 standard solution of the conjugate equations,
 $${\kappa_l(v)}^{-1}=R_v^*\otimes 1_v\circ 1_{\overline{v}}\otimes\sigma(v,v)\circ R_v\otimes 1_v,$$
 $${\kappa_l(v)}=R_v^*\otimes 1_v\circ 1_{\overline{v}}\otimes\sigma(v,v)^{-1}\circ R_v\otimes 1_v.$$
 The same relations hold for a right braided symmetry with ${\kappa_r(v)}$ in place of ${\kappa_l(v)}$.
 }\medskip
 
 \noindent{\it Proof} 
 $$R_v^*\otimes 1_v\circ 1_{\overline{v}}\otimes\sigma(v,v)\circ R_v\otimes 1_v=R_v^*\otimes 1_v\circ\sigma(\overline{v}, v)^{-1}\otimes 1_v\circ\sigma(\overline{v}\otimes v, v)\circ R_v\otimes 1_v=$$
 $$(\sigma_d(\overline{v}, v)\circ R_v)^*\otimes 1_v\circ 1_v\otimes R_v=\overline{R}_v^*\otimes 1_v\circ \kappa_l(v)^{-1}\otimes R_v=\kappa_l(v)^{-1}.$$
 The second relation follows applying the first to $\sigma_d$.
 For the remaining statements it suffices to replace $\sigma$ with $\sigma_*$ and $\sigma_{-1}$.
 \medskip
 
The following known result exhibits  the behaviour of $\kappa$ under 
tensor products, see  \cite{LR}.
 \medskip
 
 \noindent{\bf 4.4. Proposition} {\sl If $\sigma$ is a braided symmetry of ${\cal A}$, then for any pair of objects $u,z$,
 $$\kappa(u\otimes z)=(\sigma(z,u)\circ\sigma(u,z))^{-1}\circ \kappa(u)\otimes\kappa(z).$$}\medskip
 
The following convention suggested by categories of endomorphisms of an algebra helps to simplify notation. 
 The symbol $T$ denoting an arrow in $(w,z)$ will also be used for $T\otimes 1_u\in(w\otimes u, z\otimes u)$. The meaning of $T$ will be clear from the context. 
Whereas $1_u\otimes T$ will be denoted by $\rho(T)$. 
 The previous proposition easily yields the following formula.
\medskip

\noindent{\bf 4.5. Corollary} {\sl For any object $u\in{\cal A}$ and any integer $n$,
$$\kappa(u^{\otimes n})=\Sigma_{n-1}^{-1}\circ\dots\circ\Sigma_1^{-1}\circ\kappa(u)^{\otimes n},$$
where
$$\Sigma_k:=\rho^{k-1}(\sigma)\circ\rho^{k-2}(\sigma)\dots\circ\sigma^2\circ\dots\circ\rho^{k-2}(\sigma)\circ\rho^{k-1}(\sigma)$$
and $\sigma=\sigma(u,u)$.
}\medskip

\noindent{\it Remark} As $\kappa(u)$ is a scalar, when $u$ is irreducible, $\kappa(v)$ is an eigenvalue of  
$\Sigma_{n-1}^{-1}\circ\dots\circ\Sigma_1^{-1}\kappa(u)^n,$ 
 if $v$ is an irreducible summand of $u^{\otimes n}$. 
 
\medskip

 We next compare the invariant $\kappa$  of a braided tensor $C^*$--category with 
 the twist $\theta$ of a ribbon category. 
Recall that in a ${\mathbb C}$--linear tensor category (always assumed 
abelian and semisimple, although not 
necessarily strict)
a {\it right 
dual} $u^*$ 
for any object $u$ of the category,
is defined by two arrows $e\in(v^*\otimes v,\iota)$ and
$d\in(\iota, v\otimes v^*)$ such that, up to canonical 
associativity isomorphisms, omitted here for simplicity,
$$1_v\otimes e\circ d\otimes 1_v=1_v,$$
$$e\otimes 1_{v^*}\circ1_{v^*}\otimes d=1_{v^*}.$$
As for conjugates, a  right dual is unique up to a unique invertible 
$T\in(v^*_1, v^*_2)$ such that $e_1=e_2\circ T\otimes 1_v$,
$d_2=1_v\otimes T\circ d_1.$ 
A {\it left dual} $^*v$ of $v$ is similarly defined by arrows
$e'\in(v\otimes^*v,\iota)$, $d'\in(\iota, {}^*v\otimes v)$.
In a tensor $^*$--category with conjugates,
a conjugate $\overline{v}$ of $v$ is always a right and left dual, 
$\overline{v}=v^*={}^*v$.

A {\it right rigid} tensor category is a tensor category with  
a specified choice of right dual $(v^*, e_v, d_v)$ for every object $v$.
A {\it ribbon} category is a right rigid braided tensor category
with a choice  of isomorphisms $\theta_v\in(v,v)$, called {\it twists},
natural in $v$ and satisfying
$$\theta_{v\otimes 
w}=\sigma(w,v)\circ\sigma(v,w)\circ\theta_v\otimes\theta_w,$$
$$\theta_{v^*}=(\theta_v)^c,$$
$$\theta_\iota=1_\iota.$$
In a ribbon category  we also have an associated left duality $({}^*v=v^*, 
e'_v, d'_v)$
defined as in equation $(3.5)$ of Ch. XIV in  \cite{Kassel}.
The associated contravariant functor coincides with  that induced by right duality.
As is well known,  we may then define a scalar valued  trace ${\rm tr}_v$ as in Def. XIV.4.1 in \cite{Kassel}, 
analogous to a 
left inverse in a tensor $C^*$--category \cite{LR}.
The twist can be computed from the braided symmetry
and trace, 
$(\theta_v)^{-1}={\rm tr}_v\otimes 1_v(\sigma_{-1}(v,v)).$
From this, it is easy   to show that
if ${\cal N}$ is a ribbon category and
 ${\cal M}$ a tensor $C^*$--category with conjugates embedded  in ${\cal N}$ as   a full tensor subcategory then  for any irreducible object $v$ of ${\cal M}$, the trace arising from the ribbon structure may be chosen positive and
$\theta_v= \kappa_r(v)^{-1}$, cf.   Ch. XIV   in \cite{Kassel}.

\medskip

\section{An obstruction to extending  braided tensor $C^*$--categories}

Throughout  this section, ${\cal A}$ and ${\cal M}$ are tensor 
$C^*$--categories with conjugates and irreducible tensor units and 
$\rho:{\cal A}\to{\cal M}$ is a tensor $^*$--functor. 

We start by assuming  that ${\cal A}$ and ${\cal M}$ have    
braidings, $\sigma$ and $\varepsilon$, in the sense of $(2.1)$ and that $\rho$ is a braided tensor 
functor. We shall refer ${\cal M}$ as a {\it braided extension} of ${\cal A}$.

Our main assumption is that ${\cal A}$ admits an embedding into the category of Hilbert spaces, and 
we fix a tensor $^*$--functor $\tau:{\cal A}\to\text{Hilb}$.
(We shall not assume that $\tau$ is braided as this would imply that $\sigma$ is a permutation symmetry.) Hence, by Woronowicz duality, $\tau$ determines a compact quantum group $G^\tau$ with a  braided representation category.

Consider the ergodic $C^*$--action  $({\cal C},\alpha)$ of $G^\tau$ associated with the pair $(\rho,\tau)$, see \cite{PR1}. For each object $u$ of ${\cal A}$, ${\cal H}_u$ is the Hilbert bimodule constructed in \cite{PR},  in fact just depending on $\rho_u$ of ${\cal M}$.

 We identify ${\cal H}_u$ with 
$\tau_u\otimes{\cal C}$ as right Hilbert bimodules.  The canonical unitaries 
 $S_u$   make the left module structures explicit, see Prop. 8.6 in \cite{PR},
 $$<\psi\otimes I, c^v(\phi)\cdot \psi'\otimes I>=(\rho(R^*_u)\circ 1_{\rho_{\overline{u}}}\otimes 
T\otimes 1_{\rho_u})\otimes 
(j_u\psi\otimes\phi\otimes \psi'),$$ 
for every irreducible spectral representation $v$ of the ergodic action of $G^\tau$ on ${\cal C}$ and every linear
intertwiner $c^v: \tau_v\to{\cal C}$ between $v$ and $\alpha$ i.e.\ of the form
$c^v(\phi)=T\otimes\phi$, $T\in(\rho_v,\iota)$. 

The next lemma shows that if ${\cal A}$ and ${\cal M}$ have a braiding in the weak sense of relation $(2.1)$ and  if  $\rho$ is a {\it braided} functor,
the left bimodule structure is completely determined by the representation theory of $G^\tau$ and the ergodic action.

\medskip

\noindent{\bf 5.1. Lemma} {\sl If $\sigma$ and $\varepsilon$ satisfy 
$(2.1)$ for  tensor $C^*$--categories ${\cal A}$ and ${\cal M}$  and if $\rho:{\cal A}\to{\cal M}$ is a braided tensor $^*$--functor 
 then the left module structure on ${\cal H}_u$ under the canonical identification with $\tau_u\otimes{\cal C}$ is given by
$$<\psi\otimes I, c^v(\phi)\cdot \psi'\otimes I>=c^v(\tau(R_u^*\otimes 1_v\circ 1_{\overline u}\otimes \sigma(v,u)) j_u\psi\otimes\phi\otimes\psi').$$
}\medskip

\noindent{\it Proof}  
 Writing, as above,  $c^v(\phi)=T\otimes\phi$, with $\phi\in \tau_v$ and $T\in(\rho_v,\iota)$,
 $$<\psi\otimes I, c^v(\phi)\cdot \psi'\otimes I>=(\rho(R^*_u)\circ 1_{\rho_{\overline{u}}}\otimes 
T\otimes 1_{\rho_u})\otimes 
(j_u\psi\otimes\phi\otimes \psi')=$$ 
 $$(\rho(R^*_u)\circ 1_{\rho_{\overline{u}}}\otimes 
(1_{\rho_u}\otimes T\circ\varepsilon(\rho_v,\rho_u)))\otimes 
(j_u\psi\otimes\phi\otimes \psi')=$$ 
$$(\rho(R^*_u)\circ 1_{\rho_{\overline u}}\otimes 1_{\rho_u}\otimes T\circ
1_{\rho_{\overline u}}\otimes 
\rho\sigma(v,u))\otimes(j_u\psi\otimes\phi\otimes\psi')=$$
$$(T\circ\rho(R^*_u)\otimes 1_{\rho_v}\circ 1_{\rho_{\overline u}}\otimes
 \rho\sigma(v, u))\otimes(j_u\psi\otimes\phi\otimes\psi')=$$
 $$T\circ\rho(R_u^*\otimes 1_v\circ
 1_{\overline 
u}\otimes\sigma(v,u))\otimes(j_u\psi\otimes\phi\otimes\psi')=$$
$$T\otimes(\tau(R_u^*\otimes 1_v\circ 1_{\overline u}\otimes \sigma(v,u)) 
j_u\psi\otimes\phi\otimes\psi'))=$$
$$c^v(\tau(R_u^*\otimes 1_v\circ 1_{\overline u}\otimes \sigma(v,u)) j_u\psi\otimes\phi\otimes\psi').$$
 \medskip

From Lemma 5.1, we derive a first property 
that spectral representations 
for ergodic actions arising from braided functors 
$\rho:{\cal A}\to{\cal M}$ have to satisfy.

\medskip

\noindent{\bf 5.2. Corollary} {\sl 
Let $\sigma$ and $\varepsilon$ be weak left   
braided symmetries of ${\cal A}$ and 
${\cal M}$ respectively and
$\rho:{\cal A}\to{\cal M}$  a braided tensor $^*$--functor. 
For every irreducible object  $v$ of ${\cal A}$ such that 
$(\iota,\rho_v)\neq\{0\}$
and for every object
$u\in{\cal A}$, and every solution $R_u$ of the conjugate equations for $u$,
$$R_u^*\otimes 1_v\circ1_{\overline{u}}\otimes\sigma(v,u)=R_u^*\otimes 
1_v\circ1_{\overline{u}}\otimes\sigma_d(v,u).\eqno(5.1)$$
}\medskip

\noindent{\it Proof} It suffices to apply Lemma 5.1 to $\sigma$ and 
$\sigma_d$, recalling that 
entries of spectral multiplets of an ergodic action corresponding to irreducible representations are linearly independent and that $\tau$ is faithful on arrows being defined on a tensor $C^*$--category with conjugates.\medskip

\noindent{\it Remark} When $\sigma$ and $\varepsilon$ are weak right  braided symmetries,
$(5.1)$ is replaced by
$R_u^*\otimes 1_v\circ1_{\overline{u}}\otimes\sigma_{-1}(v,u)=R_u^*\otimes 
1_v\circ1_{\overline{u}}\otimes\sigma_*(v,u).$
\medskip

 Note that  $(5.1)$ is automatically satisfied if $\sigma$ is 
unitary. However, if it is non-unitary,  the corollary provides  an obstruction to 
constructing non-trivial braided extensions $({\cal M}, \rho)$ of ${\cal A}$. 
\medskip

We start with a given 
braided extension $({\cal M},\rho)$  and draw conclusions 
about the eigenvalues $\kappa_l(v)$ corresponding to irreducible objects $v$ of ${\cal 
A}$
with $(\iota, \rho_v)$ non-trivial.

\medskip

\noindent{\bf 5.3. Corollary} {\sl Let $\sigma$  be a  left 
(right) braided 
symmetry of ${\cal A}$, $\varepsilon$ a weak left (right)    braided symmetry of  ${\cal M}$ and $\rho$ a braided tensor 
$^*$--functor.
 For any   irreducible object 
 $v$ of ${\cal A}$ such that $(\iota,\rho_v)\neq\{0\}$,
$\kappa_l(v)$ ($\kappa_r(v)$) is a phase.}\medskip

\noindent{\it Proof}  Replacing $\sigma$ and $\varepsilon$  by $\sigma_{-1}$ and $\varepsilon_{-1}$ if necessary, we may assume, by Prop.\ 4.2, that $\sigma$ and 
$\varepsilon$ are  left
braided symmetries.  
 We claim that for any irreducible $v$ with $(\iota,\rho_v)\neq0$ and any  
solution $\overline{R}_v$ of the conjugate equations for $v$,
$\sigma_*(v,\overline{v})\circ 
\overline{R}_v=\sigma_{-1}(v,\overline{v})\circ 
\overline{R}_v.$
This would show that, choosing standard solutions, $\kappa_r^{\sigma_*}(v)=\kappa_r^{\sigma_{-1}}(v)$ and the first two relations in Prop.\ 4.2 imply
$|\kappa_l(v)|=1$.
 Now applying the previous proposition to  $\varepsilon$ and $\sigma$, 
 being left and hence  weak left  braided symmetries, it follows from  
$(5.1)$ 
with $v$ and $\overline{v}$ in place of $u$ and $v$ respectively,
that 
$${R}_v^*\otimes 1_{\overline v}=
{R}_v^*\otimes 1_{\overline v}\circ 1_{\overline{v}}\otimes 
\sigma_d(\overline{v}, v)\sigma(\overline{v}, v)^{-1}.$$
Composing on the right by 
$1_{\overline v}\otimes 
{\overline R}_v$ gives
$$1_{\overline v}={R}_v^*\otimes 1_{\overline v}\circ 
1_{\overline 
v}\otimes(
\sigma_d(\overline{v}, v)\circ\sigma(\overline{v}, v
)^{-1}\circ 
\overline{R}_v).$$
On the other hand
$\sigma_d(\overline{v}, v)\circ\sigma(\overline{v}, v)^{-1}\circ\overline{R}_v$
must be a scalar
multiple of $\overline{R}_v$ as $v$ is irreducible. Our equation then 
shows that the 
scalar equals $1$, so $\sigma(\overline{v}, v)^*\circ\overline{R}_v=\sigma(\overline{v}, v)^{-1}\circ\overline{R}_v$.
\medskip

\noindent{\bf 5.4. Theorem} {\sl Let $({\cal A},\sigma)$ be a  tensor 
$C^*$--category 
with conjugates and a left (right) braided symmetry. Assume that ${\cal A}$
 admits a tensor $^*$--embedding $\tau$ into the category of Hilbert spaces. 
If $\kappa(v)$ is not a phase
whenever $v$ is an irreducible of ${\cal A}$ not
equivalent to $\iota$, then every braided tensor $^*$--functor $\rho:({\cal A},\sigma)\to({\cal M},\varepsilon)$ into a weak left (right)  braided tensor $C^*$--category  $({\cal M},\varepsilon)$ is full.}\medskip

\noindent{\it Proof} The arrow space   $(\iota,\rho_v)$ of ${\cal M}$ can be identified with the spectral space of  $\hat{v}$ for the  action of $G^{\tau}$ on ${\cal C}$. We may obviously replace ${\cal M}$ by the 
category whose objects are those of ${\cal A}$ and where the space of arrows from $u$ to $v$ is now 
$(\rho_u,\rho_v)$ defining $\rho$ and the algebraic structure in the obvious way. 
By the main result of \cite{PR}, Theorem 7.7, there is then a full and 
faithful  $^*$--functor $\lambda$ from ${\cal M}$  to the category of bimodule $G^\tau$--representations, taking 
$u$ to the $G^{\tau}$--bimodule $\tau_u\otimes{\cal C}$ and
$\rho(T)$   to  $\tau(T)\otimes I.$
By the previous corollary, $(\iota,\rho_v)=\{0\}$ 
for every irreducible $v\neq\iota$, hence ${\cal C}={\mathbb C}$. Now regarding $\lambda$ as taking values 
 in the category of $G^\tau$--representations,  $\lambda\rho=\tau$ and is full, thus
$\rho$ is full. 
\bigskip

\section{Proof of Theorem  3.1}

By Theorem 5.4 it suffices to show  that $\kappa(v)$ is not a phase whenever 
$v$ is an irreducible unitary representation of $G_\mu$ not equivalent to $\iota.$

We want to compute  the invariant $\kappa$ for the braided $C^*$--tensor category $(\text{Rep}S_\mu U(d),\sigma_\omega)$ and begin with its value on  the fundamental  representation $u$, writing $\sigma$ for $\sigma_\omega$ for brevity.

Recall that $\overline{u}$ may be realized as the 
subobject of $u^{\otimes d-1}$ defined by 
$\eta(E_{d-1})$. We have the relations 
 $g_iS=-q S$,  
with $q=\mu^2$ for $i=1,\dots, d-1$
(see e.g. Lemma 5.6 in \cite{P}), hence $\sigma_i S=-\omega\mu S$. 
 
By Theorem 5.5 in \cite{P}, a standard solution of the conjugate equations 
 for the fundamental representation $u$ is given by
 $R=\lambda S$, $\overline{R}=(-1)^{d-1}R$, with $\lambda$ a suitable positive scalar.
 Hence 
 $$\sigma(u,\overline{u})\circ \overline{R}=(-1)^{d-1}\lambda\sigma(u,\overline{u})\circ 1_u\otimes E_{d-1}\circ S=$$
 $$(-1)^{d-1}\lambda E_{d-1}\otimes 1_u\circ \sigma(u, u^{\otimes d-1}) \circ S=
 (-1)^{d-1}(-\omega\mu)^{d-1}\lambda  E_{d-1}\otimes 1_u\circ S=$$
 $$(\omega\mu)^{d-1}R,$$
 so that $$\kappa(u)=(\omega\mu)^{d-1}.$$
 
 \noindent{\it Remark}
 A similar computation shows that $\kappa(\overline{u})=(\omega\mu)^{d-1}$ as well.
 \medskip
 
Note that  these are not phases, unless   $\mu=\pm 1$.
As we shall see more precisely later, these eigenvalues can be well understood in terms of 
the representation theory of the quantum group.

From the above computation of $\kappa(u)$,  
and the remark following Corollary 4.5, we arrive at the following 
result expressed 
in terms of the original representation $\eta$ 
of the Hecke algebra $H_n(\mu^2)$  

\medskip

\noindent{\bf 6.1. Proposition} {\sl
$$\kappa(u^{\otimes n})=(\frac{\omega}{\mu})^{-n(n-1)}(\omega\mu)^{n(d-1)}\eta_\mu(G_{n-1}^{-1}\dots G_1^{-1}),$$
where $G_k=g_kg_{k-1}\dots g^2\dots g_{k-1}g_k$.}\medskip

We shall use this formula to reduce the problem to the case $\mu>0$ and a specified 
root.
\medskip
 
\noindent{\bf 6.2. Theorem} {\sl  Suppose that, for some $\mu>0$ and a 
specified $d^{\text{th}}$--root of $\omega$ of $\mu$,
$\kappa(v)$ is never a phase when $v$ is an irreducible object of  $\text{Rep}(S_\mu U(d))$ not equivalent to $\iota$ then the same property holds if $\mu$ and $\omega$ are replaced independently   by  $-\mu$   and any
other root $\omega'$.}\medskip

\noindent{\it Proof} If $v$ is an irreducible summand 
of some tensor power $u^{\otimes n}$, $\kappa(v)$ appears as the eigenvalue
of  $\kappa(u^{\otimes n})$ corresponding 
to the central support of $v$ in $(u^{\otimes n}, u^{\otimes n})=\eta_\mu(H_n(q))$. 
On the other hand the representations $\eta_\mu$ and $\eta_{-\mu}$ have the same kernel $I_n$ in each $H_n(q)$ (cf. \cite{KW} and also  \cite{P}), hence this explicit dependence on the Hecke algebra $H_n(\mu^2)$ shows that the relevant  central support does not change if one replaces $\mu$ by $-\mu$.
 This central support agrees with that of an irreducible summand $v'$ in $u'^{\otimes n}$ where $u'$ is
 the fundamental representation of $S_{-\mu}U(d)$. The spectrum 
 of $\eta_\mu(G_n^{-1}\dots G_1^{-1})$ is just  the spectrum
 of the image of $G_n^{-1}\dots G_1^{-1}$ in $H_n(q)/I_n$ under the quotient map, and hence 
 does not change if $\mu$ is replaced by $-\mu$. It follows that
 $\kappa(v)$   for $\mu$ and $\sigma_\omega$  can differ from  $\kappa(v')$ for $ \mu$ or  $-\mu$ and $\sigma_{\omega'}$   only by a phase. \medskip

\noindent{\it Remark} The kernel of $\eta_\mu$ is known and 
may be described in terms of Young diagrams. Hence the 
spectral analysis of the element $G_{n-1}^{-1}\dots G_1^{-1}$ of the Hecke algebra $H_n(q)$ determines the invariant $\kappa$ on all irreducibles
contained in $u^{\otimes n}$, cf. \cite{WenzlInv}. However,
we shall not pursue  this approach, but rather look for a 
more general argument giving results for deformations of classical compact Lie groups of Lie types other than  $A$.
\medskip

 To this end, we use the dual picture of Drinfeld and Jimbo. See Sect.\ 3 for notation.
   
Let $\rho\in h^*$ be the element defined by $(\alpha_i, 2\rho)=(\alpha_i,\alpha_i).$
Let $\lambda$ be a dominant integral weight and $v_\lambda$ the irreducible representation of $U_\mu(g)$ with highest weight $\lambda$. The twist $\theta_{v_\lambda}$ is known to act on the space of $v_\lambda$ as scalar multiplication by $\mu^{-(\lambda,\lambda+2\rho)}$
 see Lemma 7.3.2 in
\cite{KS}.

Taking the comparison between the algebraic and analytic approach to the twist into account, see Prop.\ 5.1, we are reduced to showing that $(\lambda,\lambda+2\rho)$ is strictly positive for 
$\lambda\neq0.$ Now by \cite{Hum} Sect.\ 13.3, if
 $(\lambda_j)$ are the fundamental dominant weights  (the basis dual to $(\frac{2\alpha_j}{(\alpha_j,\alpha_j)})$),
 $\rho=\sum_1^r\lambda_j$. Writing $\lambda$ in the form $\lambda=\sum_i m_i\lambda_i$,
with $m_i$ non-negative integers, 
$(\lambda,\lambda+2\rho)=\sum_{i,j}m_i(m_j+2)(\lambda_i,\lambda_j)$.
The explicit relationship between simple roots and fundamental weights,
$\lambda_i=\sum_r d_{ir}\alpha_r$, with $(d_{ir})=({A^t})^{-1}$,
  shows that
$$(\lambda_i,\lambda_j)=\sum_{r,s}d_{ir}d_{js}(\alpha_r,\alpha_s)=\sum_{r.s}d_{ir}d_{js}d_ra_{rs}=d_jd_{ij}>0,$$
(cf. table 1, p.69, \cite{Hum}), completing the proof of Theorem 3.1.
\medskip

In particular,
 in the type $A_{d-1}$ case
that table in Humphrey's book gives,
$$(\lambda_i, \lambda_i+2\rho)=\frac{2}{d}(d-i)(1+\dots+i-1)+\frac{3}{d}i(d-i)+\frac{2}{d}i(1+\dots+d-i-1)=$$
$$\frac{1}{d}(d-i)(i-1)i+\frac{3}{d}i(d-i)+\frac{1}{d}i(d-i-1)(d-i)=$$
$$\frac{1}{d}i(d-i)(d+1),$$ 
for the fundamental weights. For $i=1$ and $i=d-1$, this 
reduces to our previous  computation of $\kappa(u)$ and $\kappa(\overline{u})$.\bigskip

\noindent{\it Remark} The computation of $\kappa(u)$ and $\kappa(\overline{u})$   at the beginning of the section depended on  relations $(3.1)$--$(3.3)$. The above proof gives an independent derivation. 
\medskip

\bigskip

\section{ On the characterization of $\text{Rep}(S_\mu U(d))$}

An intrinsic characterization of the quantum deformation of the spatial 
linear group was first given in \cite{KW}, and based on the 
analysis of the fusion rules. An independent approach was proposed in 
\cite{P} for $S_\mu U(d)$, where the starting point was the braided 
symmetry. The aim of this  section is to make the relation between  the two 
approaches  explicit, when not at roots of unity. 

To any semisimple rigid  tensor category ${\cal C}$ with Grothendieck semiring isomorphic to that of $sl(d)$ (briefly, an $sl(d)$--category), Kazhdan and Wenzl associate an invariant $\tau_{\cal C}$, the {\it  twist of the category}. They start with a suitable
 idempotent $a\in(X^2,X^2)$,  $X$ the fundamental object of ${\cal C}$, and find a nonzero
complex number $q$ such that $qa-(I-a)$ satisfies the Hecke algebra relations for the parameter $q$.
This defines  representations $\eta_{KW}:H_n(q)\to(X^n,X^n)$. If $\nu\in(\iota, X^d)$
and $p\in(X^d,\iota)$ are chosen so that $p\circ \nu=1$ then the twist  is defined by
$$\tau_{\cal C}:=p\otimes  1\circ \eta_{KW}(g_d\dots g_1)\circ 1\otimes\nu.$$
It is asserted in Prop.\ 5.1, \cite{KW}  that 
$\tau_{\cal C}$ is a $d^{\text{th}}$-root of unity.
Clearly, $\text{Rep}(S_\mu U(d))$  is an $sl(d)$--category. 
However, comparing with  $(3.3)$ leads to  $\tau=\mu^{d-1}$.
The origin of the inconsistency is
 the  claim on page 135 of \cite{KW} that $\nu\otimes\nu$ is an eigenvector of
 $\eta_{KW}(X^d, X^d)$ with eigenvalue $1$, whilst, for $S_\mu U(d)$, iterating relations $(3.3)$   $d$ times 
 gives, 
$\eta(u^d, u^d)\nu\otimes \nu=\mu^{d(d-1)}\nu\otimes \nu.$

This value is related to Drinfeld's theory of universal  $R$--matrices. 
In fact, it is well known that universal
 $R$--matrices give rise to ribbon tensor categories. As a consequence, in 
the type $A$ case,
the generator of the Hecke algebra needs to be multiplied by a 
suitable scalar to ensure the naturality of the braiding. 
This scalar is 
well known, see Lemma 3.2.1 in \cite{WenzlInv}, see also 
\cite{Blanchet}, \cite{Y}.
(We would warn the reader that this scalar is often computed incorrectly  
in the literature.)

The well known characterization of  the 
quasiequivalence class of a non-faithful Hecke algebra representation in 
a rigid tensor 
category ${\cal C}$ \cite{KW, W, P} leads to $\tau_{\cal 
C}=w\mu^{d-1}$ where $\mu$ 
is a 
complex square root of $q$ and $w$ is a $d^{\text{th}}$  root of unity.

This expression easily leads to a
presentation of $sl(d)$--categories analogous to \cite{P}.
To make the comparison more immediate, we shall assume 
that ${\cal C}$ has a ${}^*$--involution making it into a 
tensor ${}^*$--category. Analogous  relations may be derived in the general case.
We omit the proof.
\medskip

\noindent{\bf 7.2. Proposition} {\sl If ${\cal C}$ is  a tensor 
${}^*$--category of 
$sl(d)$ type,   $q$ is derived from ${\cal C}$ as in \cite{KW}, 
$\nu\in(\iota, X^d)$ satisfies
$\nu^*\circ\nu=1$,
then 
$ \nu\circ\nu^*=\eta(E_d)$ and
$$\nu^*\otimes 1_X\circ 1_X\otimes 
\nu=w\frac{(-\mu)^{d-1}}{[d]_q},\eqno(7.2)$$
$$\eta(g_1\dots g_d)\nu\otimes 1_X=\overline{w} \mu^{d-1} 1_X\otimes 
\nu,\eqno(7.3)$$
where $\mu^2=q$, $\tau_{\cal C}=w\mu^{d-1}$, $[d]_q=1+q+\dots+q^{d-1}$,  and $\eta=\eta_{KW}$.
}\medskip

\noindent{\it Remark} If ${\cal C}$ is a tensor $C^*$--category, there are Hilbert space representations for $H_\infty(q)$. Hence $q$ is either a root of unity or $q>0$
by Wenzl's result \cite{Wenzl}. 
\medskip

We characterize $\text{Rep}(S_\mu U(d))$ among $sl(d)$--categories in the spirit 
of \cite{KW}.
\medskip
 
 \noindent{\bf 7.3. Corollary} {\sl Let ${\cal C}$ be a  tensor 
$C^*$--category  of $sl(d)$--type with associated parameter $q$. Then,
  \begin{description}
  \item{\rm a)}
  if $\tau_{\cal C}>0$ then ${\cal C}$ is tensor $^*$--equivalent to $\text{Rep}(S_{\sqrt{q}} U(d))$,
\item{\rm b)} if $\tau_{\cal C}<0$ and $d$ is even then ${\cal C}$ is tensor 
$^*$--equivalent to
$\text{Rep}(S_{-\sqrt{q}} U(d))$.
\end{description}

}\medskip

\noindent{\it Remark} Note that for $\mu>0$, and $d$ even, 
$\text{Rep}(S_{-\mu}U(d))$ is a twist of $\text{Rep}(S_\mu U(d))$ by
$w=-1.$
\medskip

 \noindent{\bf Acknowledgements.} C.P. would like to thank  C. De 
Concini, Y. Kawahigashi, R. Longo  and M. M\"uger for discussions. We would like to thank the 
referee for    pointing out Lemma 3.2.1 in \cite{WenzlInv} and \cite{Blanchet} and for giving suggestions on how  to shorten the presentation of the paper.
 
\bigskip 

\noindent  Dipartimento di Matematica, Universit\`a di Roma ``La Sapienza'',
00185--Roma, Italy, pinzari@mat.uniroma1.it
\bigskip

\noindent  Dipartimento di Matematica, Universit\`a di Roma ``Tor Vergata'',
00133--Roma, Italy, roberts@mat.uniroma2.it
\bigskip

\bigskip

\begin{thebibliography} {VD}

\bibitem{A} N. Andruskiewitsch: Some exceptional compact matrix pseudogroups,
{\it Bull. Soc. Math. France\/}, {\bf 120} (1992), 297--325.


\bibitem{AH} M. Asaeda, U. Haagerup: Exotic subfactors of finite depth with Jones indices 
$(5+\sqrt{13})/2$ and $(5+\sqrt{17})/2$, {\it  Comm. Math. Phys.\/}, {\bf 202}   (1999),   1--63.


\bibitem{BK} B. Bakalov, A. Kirillov, Jr.: Lectures on tensor categories and modular functors. 
University Lecture Series, vol. 21, AMS 2001.

\bibitem{B} T. Banica: Le groupe quantique compact libre $U(n)$, {\it Commun. Math. Phys.\/}, {\bf 190} (1997), 143--172.

\bibitem{Blanchet} C. Blanchet: Hecke algebras, modular categories and 
$3$-manifolds quantum invariants, {\it  Topology},  {\bf 39}  (2000),  
193-223. 

\bibitem{CKL} S. Carpi, Y. Kawahigashi, R. Longo: On the Jones index values for conformal subnets,
arXiv:1002.3710.

\bibitem{CP} V. Chari, A. Pressley: A guide to quantum groups, Cambridge university press, 1994.

\bibitem{D} P. Deligne: Cat\'egories tannakiennes, The Grothendieck Festschrift, Vol. II,  111--195, Progr. Math., 87, Birkh???user Boston, Boston, MA, 1990.

\bibitem{DHR} S.\ Doplicher, R. Haag, J.E.\ Roberts: Local observables and particle statistics I, {\it Commun. Math. Phys.}, {23} (1971), 199--230; II, {\it Commun. Math. Phys.}, {\bf 35} (1974), 49--85.

\bibitem{DR} S.\ Doplicher, J.E.\ Roberts: A new duality theory for compact groups, {\it Invent.\ Math.}, {\bf 98} (1989), 157--218.

\bibitem{ENO} P. Etingof, D. Nikshych, V.  Ostrik: On 
fusion categories, {\it  Ann. of Math.}, {\bf  162}  (2005),  
581--642. 

\bibitem{FRS} K. Fredenhagen, K.H. Rehren, B. Schroer: Superselection sectors with braid group statistics and exchange algebras 1. General theory, {\it Commun. Math. Phys.\/}, {\bf 125} (1989), 201--226.

\bibitem{GHJ}   F. M. Goodman, P. de la Harpe, V. F. R. Jones: 
 Coxeter graphs and towers of algebras. MSRI Publications {\bf 14}, Springer-Verlag, Berlin and New York, 1989.



\bibitem{GW} F.M. Goodman, H. Wenzl: Ideals in the Temperley-Lieb category. Ap- 
pendix to M. Freedman: a magnetic model with a possible Chern-Simons 
phase. {\it Comm. Math. Phys.\/}, {\bf  234} (2003), 129--183. 


\bibitem{Haag} R. Haag: Local quantum physics, Springer--Verlag, 1992.

\bibitem{Hum} J.E. Humphreys: Introduction to Lie algebras and representation theory.
Springer--Verlag, 1972.


\bibitem{H} U. Haagerup: Principal graphs of subfactors in the index range $4<[M:N]<3+\sqrt{2}$.
On: Subfactors (Kyuzeso, 1993), World Sci. (1994), 1--38.

\bibitem{J} M. Jimbo: A $q$--analogue of $U(gl(N+1))$. Hecke algebras and the Young Baxter equation, {\it Lett. Math. Phys.\/}, {\bf 10} (1985), 63--69.



\bibitem{Jones} V.F.R. Jones: Index for subfactors, {\it Invent. Math.\/}, {\bf 72} (1983),
1--25.

\bibitem{Kassel} C. Kassel: Quantum groups, Springer--Verlag, (1995). 


\bibitem{KRT} C. Kassel, M. Rosso, V. Turaev: Quantum groups and knot invariants,
Soc. Math. de France, 1997.

\bibitem{KL} Y. Kawahigashi, R. Longo: Classification of local conformal nets. Case $c<1$, {\it Ann. 
Math.\/}, {\bf 160} (2004), 493--522.



\bibitem{KW} D. Kazhdan, H. Wenzl: Reconstructing monoidal categories. I. M. Gelfand seminar,
{\it Adv. Soviet. math.\/}, {\bf 16}, AMS, Providence, RI, (1993) 111--136.

\bibitem{KS} L.I. Korogodski, Y.S. Soibelman: Algebras of functions on quantum groups: part I, Mathematical Surveys and Monographs, AMS, {\bf 56}, 1998.


\bibitem{LR} R.\ Longo, J.E.\ Roberts: A theory 
of dimension, {\it
$K$--Theory}, {\bf 11} (1997), 103--159.



\bibitem{Muger} M. M\"uger: On the structure of modular categories, {\it Proc. London Math. Soc.\/}
{\bf 87} (2003), 291--308.

\bibitem{muger} M. M\"uger: On superselection theory of quantum fields in low dimensions, arXiv:0909.2537.

\bibitem{Mlectures} M. M\"uger: Tensor categories: A selective guided tour, arXiv:0804.3587.

\bibitem{O1} V. Ostrik: Fusion categories of rank 2, {\it Math. Res. Lett\/}, {\bf 10} (2003), 177--183.

\bibitem{O2} V. Ostrik: Pre-modular categories of rank 3, {\it Mosc. Math. J.\/}, {\bf 8} (2008), 
111--118.

\bibitem{P} C. Pinzari: The representation category of the Woronowicz quantum group $S_\mu U(d)$ as a braided tensor $C^*$--category, {\it Int. J. Math.\/}, {\bf 18} (2007), 113--136.

\bibitem{PR1} C.\ Pinzari, J.E.\ Roberts: A duality theorem for ergodic
actions of compact quantum groups on $C^*$--algebras, {\it Comm.\ Math.\ 
Phys.}, {\bf 277} (2008). 385-421.


\bibitem{PR} C. Pinzari, J.E. Roberts: A theory of induction and classification of tensor $C^*$--categories, arXiv:0907.2459.

\bibitem{Popa} S. Popa: Classification of amenable subfactors of type $II$,
{\it Acta Math.}, {\bf 172} (1994), 163--255.

\bibitem{Rehren} H. Rehren: Subfactors and coset models. In: Generalized symmetries in Physics, 
Clausthal 
(1994), 338--356.

\bibitem{RT} N. Yu. Reshetikhin, V.G. Turaev:  Ribbon graphs and their invariants derived from quantum groups, {\it Comm. Math. Phys.}, {\bf 127} (1990), 1--26.

\bibitem{R} M. Rosso: Alg\`ebres enveloppantes quantifi\'ees, groupes quantiques compacts de matrices et calcul diff\'erentiel non commutatif, {\it Duke Math. J. \/},
{\bf 61}  (1990), 11--40.


\bibitem{RSW} E. Rowell, R. Stong, Z. Wang: On classification of modular tensor categories, {\it Commun. 
Math. Phys.\/} {\bf 292} (2009), 343--389.


\bibitem{TW} I. Tuba, H. Wenzl: On braided tensor categories of type BCD, {\it J. Reine Angew. Math.\/},
{\bf 581} (2005), 31--69.

 \bibitem{Wenzl} H. Wenzl: Hecke algebras of type $A_n$ and subfactors,  {\it Invent. Math.},  {\bf 92}  (1988),  349--383.
 
 \bibitem{WenzlInv} H. Wenzl: Braids and invariants of $3$--manifolds, {\it Invent. Math.},
 {\bf 114} (1993), 235--275.

\bibitem{W} S. L.\ Woronowicz: Tannaka--Krein duality for compact matrix pseudogroups. Twisted ${\rm SU}(N)$ groups, {\it  Invent.\ Math.}, {\bf  93}  (1988),   35--76. 

\bibitem{Y} S. Yamagami: Free products of semisimple tensor categories, 
arXiv:math/0106214.

\end{thebibliography}
\end{document}